\documentclass[11pt]{article}
\usepackage{color}
\usepackage{amsmath,amssymb,comment}
\usepackage{amsfonts}
\newtheorem{theo}{Theorem}[section]
\newtheorem{lem}[theo]{Lemma}

\newcommand{\proof}{{\sc Proof.} \quad}
\newcommand{\R}{\mathbb{R}}

\newcommand{\be}{\begin{equation} \label}
\newcommand{\ee}{\end{equation}}
\newcommand{\bea}{\begin{eqnarray}\label}
\newcommand{\eea}{\end{eqnarray}}
\newcommand{\bas}{\begin{eqnarray*}}
\newcommand{\eas}{\end{eqnarray*}}
\newcommand{\bit}{\begin{itemize}}
\newcommand{\eit}{\end{itemize}}
\newcommand{\nn}{\nonumber}

\newcommand{\abs}{\\[3mm]}

\newcommand{\mysection}[1]{\section{#1} \setcounter{equation}{0}}
\newcommand{\parab}{{\cal{P}}}
\newcommand{\qarab}{{\cal{Q}}}
\newcommand{\chip}{\chi^{(\xi_0,t_0,A)}}
\newcommand{\fp}{f^{(\xi_0, t_0)}}
\newcommand{\phip}{\varphi^{(\xi_0, t_0, A)}}
\newcommand{\ophia}{\overline{\varphi}^{(A)}}
\newcommand{\hphi}{\widehat \Phi}
\newcommand{\hchia}{\widehat{\chi}^{(a)}}
\newcommand{\uphia}{\underline{\varphi}^{(a)}}

\def\qed{\hfill
\raise -2pt
\hbox{\vrule \vbox to8pt{\hrule width 6pt
\vfill\hrule}\vrule}\\[5mm]}

\begin{document}

\title{ Rate of Convergence to Barenblatt Profiles \\ for the Fast
Diffusion Equation \\
with a Critical Exponent
}


\author{Marek Fila \\
Department of Applied Mathematics and Statistics, Comenius
University \\
84248 Bratislava, Slovakia \\
\\
J.R. King   \\
Division of Theoretical Mechanics, University of Nottingham\\
Nottingham NG7 2RD, UK\\
\\
and \\
\\
Michael Winkler \\
Institut f\"ur Mathematik, Universit\"at Paderborn \\
33098 Paderborn, Germany \\
}

\date{}

\maketitle

\begin{abstract}
 We study the asymptotic behaviour near extinction of positive
solutions of the Cauchy problem
for the fast diffusion equation with a critical exponent.
After a suitable rescaling which yields a non--linear Fokker--Planck equation, we
find a continuum of algebraic rates of convergence to a self--similar profile.
These rates depend explicitly on the spatial decay rates of initial data. This improves a previous result
on slow convergence for the critical fast diffusion equation ({\sc Bonforte et al}. in Arch Rat Mech Anal
196:631--680, 2010) and provides answers to some open problems.
\end{abstract}


\mysection{Introduction}
\label{intro}
We consider the Cauchy problem for the fast diffusion equation,
\begin{equation}\label{0}
       \left\{ \begin{array}{ll}
       u_\tau = \nabla\cdot (u^{m-1}\nabla u),
               \qquad & y\in\R^n, \ \tau\in (0,T), \\[2mm]
       u(y,0)=u_0(y)\ge 0, \qquad & y\in\R^n,
       \end{array} \right.
\end{equation}
where $n\ge 3$, $T>0$ and $m=(n-4)/(n-2)$.
It is known that for $m<m_c:=(n-2)/n$ all
solutions with initial data in some suitable space, such as $L^p(\R^n)$ with
$p=n(1-m)/2$, extinguish in finite time. We shall
consider solutions which vanish in a finite time $\tau=T$
and study their behaviour near $\tau=T$.

For the extinction range $m<m_c$ there are (infinite--mass) solutions of
the self--similar form
\begin{equation}\label{baren1}
U_{D, T}(y,\tau):=\frac 1{R(\tau)^n} \left(D+\frac{\beta(1-m)}{2}
\left|\frac{y}{R(\tau)}\right|^2\right)^{-\frac{1}{1-m}},
\end{equation}
where  $D\ge 0$ and
\[
R(\tau):=(T-\tau)^{-\beta},\qquad
\beta:= \frac1{n(1-m)-2}=\frac{1}{n\,(m_c-m)}>0.
\]
 We will call these  solutions  \emph{Barenblatt solutions}.

Many papers (\cite{BBDGV,BDGV,BGV,FVW,FVWY}, for example)
are concerned with the convergence of solutions
of (\ref{0}) to the Barenblatt solutions as $\tau\to T$.
More precisely, the decay rates of
\[
R(\tau)^{n}(u(\tau,y)-U_{D, T}(y,\tau))
\]
as $\tau\to T$ are discussed there. 

The  critical exponent
\[
m_*:=\frac{n-4}{n-2}<m_c\, ,
\]
has the property that the difference of two Barenblatt
solutions is integrable for $m\in (m_*,m_c)$, while it is not integrable for
$m\leq m_*$. The exponent $m_*$ plays a very important role in the
results of \cite{BBDGV,BDGV,BGV,FVW,FVWY}. 

To study the asymptotic profile as $\tau\to T$,
it is convenient to rewrite
(\ref{0}) in self--similar variables:
\bas
t:=\frac{1}{\mu}\log\left(\frac{R(\tau)}{R(0)}\right)\quad\mbox{and}\quad
x:=\sqrt{\frac{\beta}{\mu}}\,\frac y{R(\tau)}\,,\quad \mu:=\frac{2}{1-m},
\eas
with $R$ as above, and the rescaled function
\bas
v(x,t):= R(\tau)^{n}\,u(y,\tau)
\eas
satisfies then the non--linear Fokker--Planck equation
\begin{equation}\label{fp}
	v_t= \nabla\cdot (v^{m-1}\nabla v)+\mu \,\nabla\cdot(x\,v), \qquad x\in \R^n, \, t>0.
\end{equation}

The Barenblatt solutions $U_{D,T}
(y,\tau)$ are thereby transformed into \emph{Barenblatt profiles} $V_D(x)$, which
have the advantage of being stationary:
\be{vd}
V_D(x):=(D+|x|^2)^{-1/(1-m)},\quad x\in \R^n\,.
\ee
In the new variables, the convergence of solutions of (\ref{0}) to $U_{D,T}$ takes the form 
of stabilization of solutions of (\ref{fp}) to non--trivial equilibria 
$V_D$.

In \cite{BBDGV,BDGV,FVW,FVWY} one can find several sufficient conditions under which $v(\cdot,t)$
converges to $v_D$ exponentially if $m<m_c$, $m\not=m_*$. The case $m=m_*$ was treated in \cite{BGV}
by functional analytic methods. A suitable linearization of the non--linear Fokker--Planck equation (\ref{fp})
was viewed as the plain heat flow on a suitable Riemannian manifold and then non--linear stability was studied
by entropy methods. One of the main results of \cite{BGV} says that if
$0<D_1<D_0$, $D\in [D_1, D_0]$ and
\bas
V_{D_0}(x)\le v_0(x)\le V_{D_1}(x),\qquad x\in \R^n,
\eas
\be{BGV:H2}
|v_0(x)-V_D(x)|\le f(|x|),\qquad x\in \R^n,\qquad f(|\cdot|)\in L^1(\R^n),
\ee
then for the solution $v$ of (\ref{fp}) with the initial condition $v(x,0)=v_0(x)$ it holds that
\be{BGV:res1}
\Vert v(\cdot,t)-V_D\Vert_{L^\infty(\R^n)}\le K(t+1)^{-\frac{1}{4}},\qquad t\ge 0,
\ee
for some $K>0$.

No lower bound for the rate was given in \cite{BGV} and the question 
of whether the rate from
(\ref{BGV:res1}) is optimal for a class of data was posed there as 
an open problem together with the question of whether one can prove
convergence, maybe with worse rates or without rates, for more general initial data.
Our aim is to provide some answers to these questions by establishing optimal results on rates of convergence
for a class of initial data which do not satisfy (\ref{BGV:H2}).

\begin{theo}\label{thm1}
Assume that $n>2$, $m=m_*=\frac{n-4}{n-2}$ and $D>0$, and that $V_D$ is as defined in {\rm (\ref{vd})}.
Let
$v$ be the solution of {\rm (\ref{fp})} with the initial condition
\be{v_0}
	v(x,0)=v_0(x):=\Big(|x|^2+D+\psi_0(x) \Big)^{-\frac{n-2}{2}}, \qquad x\in\R^n,
\ee
where $\psi_0$ is continuous and nonnegative on $\R^n$, $\psi_0\not\equiv 0$. 
\item{\rm(i)} If there are $B>0$ and $\gamma \in (0,1)$ such that
 \bas
	\psi_0(x) \le B \ln^{-\gamma} |x|,
	\qquad |x|>2,
  \eas
  then there exists $C>0$ such that 
  \bas
	  V_D(x) \left(1- C V_D^{\frac{2}{n-2}}(x) (t+1)^{-\frac{\gamma}{2}}\right)\le v(x,t)\le V_D(x),\qquad x\in\R^n,\quad 
	 t\ge 0.
  \eas
\item{\rm(ii)} If there are $b>0$ and $\gamma \in (0,1)$ such that
  \bas
	\psi_0(x) \ge b \ln^{-\gamma}|x|,
	\qquad |x|>2,
  \eas
then there exists $c>0$ such that 
  \bas
	v(0,t) \le V_D(0) - c (t+1)^{-\frac{\gamma}{2}},
	\qquad t>0.
  \eas
\end{theo}
This theorem says that if $V_D(x)-v_0(x)$ behaves like $|x|^{-n}\ln^{-\gamma}|x|$ for $|x|$ large and some
$\gamma\in(0,1)$ then 
$\Vert v(\cdot,t)-V_D\Vert_{L^\infty(\R^n)}$ behaves like $t^{-\gamma/2}$ for $t$ large. Hence, we obtain a continuum of
algebraic rates for initial data which do not satisfy (\ref{BGV:H2}).
We shall also show that convergence to $V_D$ from below
cannot occur at any rate faster than $t^{-1/2}$, so
Theorem~\ref{thm1} (i) does not hold for $\gamma>1$.
\begin{theo}\label{thm2}
  Let $n>2, m=m_\star$ and $D>0$,
  and assume that $\psi_0$ is continuous and nonnegative on $\R^n$, $\psi_0 \not\equiv 0$.
  Then there exists $c>0$ such that the solution $v$ of {\rm (\ref{fp}), (\ref{v_0})} satisfies
  \bas
	v(0,t) \le V_D(0) - c (t+1)^{-\frac{1}{2}}
	\qquad \mbox{for all } t>0,
  \eas
 with $V_D$ as given by {\rm (\ref{vd}).}
\end{theo}

We prove Theorems \ref{thm1} and \ref{thm2} by constructing suitable sub-- and super--solutions. In order not to 
make the paper unnecessarily long, we consider only initial data below $V_D$ but one can modify the arguments to
prove analogous results for initial data above $V_D$.

In Section~\ref{Sec2}
we establish the lower bound from Theorem~\ref{thm1}~(i) and in Section~\ref{Sec3} the upper bound from 
 Theorem~\ref{thm1}~(ii). Section~\ref{Sec4} is devoted to the proof of Theorem~\ref{thm2}.
\mysection{Lower bound. Proof of Theorem \ref{thm1} (i)}\label{Sec2}
To construct a suitable super--solution we shall use the following:
\begin{lem}\label{lem1}
  Let $\gamma\in (0,1)$. Then the solution of the problem
  \be{1.1}
	\left\{ \begin{array}{l}
	\Phi''(z)+\frac{z}{2} \Phi'(z)+\frac{\gamma}{2} \Phi(z)=0, \qquad z>0, \\[1mm]
	\Phi(0)=1, \quad \Phi'(0)=0,
	\end{array} \right.
  \ee
  is positive and decreasing on $[0,\infty)$, and there exist $c>0$ and $C>0$ such that
  \be{1.2}
	cz^{-\gamma} \le \Phi(z) \le C z^{-\gamma} \qquad \mbox{for all } z\ge 1
  \ee
  and
  \be{1.3}
	-Cz^{-\gamma-1} \le \Phi'(z) \le -cz^{-\gamma-1} \qquad \mbox{for all } z\ge 1
  \ee
  as well as
  \be{1.4}
	|\Phi''(z)| \le Cz^{-\gamma-2} \qquad \mbox{for all } z\ge 1.
  \ee
\end{lem}
\proof
The solution $\Phi$ of (\ref{1.1}) can be written explicitly in the form
\bas
\Phi(z)=e^{-\zeta}{\cal M}\left(\frac{1-\gamma}{2},\frac 12,\zeta\right),\qquad \zeta:=\frac{z^2}{4},
\eas
where $\cal M$ is Kummer's function (see \cite{AS})
\bas
{\cal M}(a,b,\zeta):=1+\frac ab \zeta+\dots+\frac{a(a+1)\dots (a+k)}{b(b+1)\dots (b+k)k!}\zeta^k+\dots
\eas
and 
\be{Kum}
\zeta^{b-a}e^{-\zeta}{\cal M}(a,b,\zeta)\to\frac{\Gamma(b)}{\Gamma(a)}\quad\mbox{as}\quad\zeta\to\infty,
\ee
which yields (\ref{1.2}).

If we now rewrite the equation in (\ref{1.1}) as
\bas
\Phi''(z)+\frac{1}{2}z^{1-\gamma}\big(z^\gamma\Phi(z)\big)'=0
\eas
and use the identity
\bas
\zeta\frac{\rm d}{{\rm d}\zeta}\big(\zeta^{b-a}e^{-\zeta}{\cal M}(a,b,\zeta)\big)=
(b-a)\zeta^{b-a}e^{-\zeta}{\cal M}(a-1,b,\zeta)
\eas
(see \cite{AS}) together with (\ref{Kum}) then we obtain that
\bas
\left|\frac 12 z^{1-\gamma}\big(z^\gamma\Phi(z)\big)'\right|\leq C  z^{-\gamma-2},\qquad z\geq 1,
\eas
which implies (\ref{1.4}).

Since $\Phi$ cannot have any local minimum, one can see that $\Phi'$ is negative and  (\ref{1.3}) follows from (\ref{1.4}).
\qed

For $m=m_*$ and radial solutions $v=v(r,t)$, (\ref{fp}) becomes
\bas
	v_t=  (v^{-\frac{2}{n-2}}v_r)_{r} + \frac{n-1}{r} v^{-\frac{2}{n-2}}v_r  +(n-2)( r v_r +  n v),
	\qquad r>0, \, t>0.
\eas
If we further transform $v$ via
\bas
	v(r,t)=\Big(r^2+D+\varphi(r,t) \Big)^{-\frac{n-2}{2}}, \qquad r \ge 0, \, t \ge 0,
\eas
then after some computation it can be checked that $\varphi$ satisfies for $r>0$ and $t>0$ the equation
\be{0phi}
	\parab \varphi := \varphi_t - (r^2+D+\varphi) \Big(\varphi_{rr} + \frac{n-1}{r} \varphi_r \Big) 
	+(n-2) r \varphi_r + \frac{n-2}{2} \varphi_r^2 =0.
\ee
The change of variables
\bas
	\chi(\xi,t):=\varphi(r,t), \quad \xi:=\ln r, \qquad r>0, \, t\ge 0,
\eas
yields that
\be{0chi}
	\qarab \chi := \chi_t - \chi_{\xi\xi} 
	- e^{-2\xi} \Big\{ (D+\chi) [\chi_{\xi\xi} + (n-2)\chi_\xi] - \frac{n-2}{2} \chi_\xi^2 \Big\} =0
\ee
for  $\xi\in \R$ and $t>0$.

In a region where $r$ is appropriately large, we shall use functions of the form
\be{2.1}
	\chip (\xi,t)
	:=A(t+t_0)^{-\frac{\gamma}{2}}  \Phi \Big((\xi+\xi_0)(t+t_0)^{-\frac{1}{2}}\Big),
	\qquad \xi \ge 0, \, t\ge 0,
\ee
as (upper) comparison functions. For clarity of notation, we consider $\xi_0>0$, $t_0\ge 1$ and $A>0$
as free parameters here. We shall fix $\xi_0$, $t_0$ in Lemma~\ref{lem7} and $A>0$ in the proof of
Lemma~\ref{lem8}.
\begin{lem}\label{lem2}
  Let $\gamma \in (0,1)$. For $t_0 \ge 1, \xi_0\in\R$ and $A>0$, the function $\chi=\chip$ defined in {\rm (\ref{2.1})}
  satisfies
  \be{2.2}
	\chi_t=\chi_{\xi\xi} \qquad \mbox{for } \xi>0 \mbox{ and } t>0.
  \ee
  Moreover, there exists $t_\star>1$ with the property that, whenever $t_0>t_\star$, for any choice of $\xi_0>0$ and $A>0$
  we have
  \be{2.3}
	\chi_{\xi\xi} + (n-2)\chi_\xi \le 0
	\qquad \mbox{for all $\xi>0$ and } t>0.
  \ee
\end{lem}
\proof
  Since
  \be{2.4}
	\chi_\xi=A(t+t_0)^{-\frac{\gamma}{2}-\frac{1}{2}} \Phi'(z), \qquad 
	\chi_{\xi\xi}=A(t+t_0)^{-\frac{\gamma}{2}-1} \Phi''(z)
  \ee
  and 
  \bas
	\chi_t=-\frac{1}{2} A (t+t_0)^{-\frac{\gamma}{2}-1}  z\Phi'(z)
	-\frac{\gamma}{2}A(t+t_0)^{-\frac{\gamma}{2}-1} \Phi(z)
  \eas
  with $z:=(\xi+\xi_0)(t+t_0)^{-1/2}$, the identity (\ref{2.2}) is immediate from (\ref{1.1}).

  To verify (\ref{2.3}), we observe that, since $\Phi''(0)<0$ by (\ref{1.1}), there exists $z_0>0$ such that
  \be{2.5}
	\Phi''(z) \le 0 \qquad \mbox{for all } z\in [0,z_0].
  \ee
  Then (\ref{1.3}) and (\ref{1.4}) ensure that with some $c_1>0$ and $c_2>0$ we have
  \be{2.6}
	\Phi'(z) \le -c_1 z^{-\gamma-1} 
	\qquad \mbox{for all } z>z_0
  \ee
  and
  \be{2.7}
	\Phi''(z) \le c_2 z^{-\gamma-2} 
	\qquad \mbox{for all } z>z_0.
  \ee
  We now let $t_\star>1$ be large enough such that
  \be{2.8}
	t_\star \ge \Big( \frac{c_2}{(n-2)c_1 z_0} \Big)^2 ,
  \ee
  and claim that (\ref{2.3}) holds whenever $t_0>t_\star$, $\xi_0>0$ and $A>0$.
  Indeed, recalling (\ref{2.4}), (\ref{2.5}) and the monotonicity of $\Phi$, we easily see that in the region
  where $z=(\xi+\xi_0)(t+t_0)^{-1/2} \le z_0$, both $\chi_{\xi\xi}$ and $\chi_\xi$ are nonpositive and hence
  clearly $\chi_{\xi\xi}+(n-2)\chi_\xi \le 0$. 
  On the other hand, if $z> z_0$ then from (\ref{2.4}), (\ref{2.6}) and (\ref{2.7}) it follows that
  \bas
	\frac{\chi_{\xi\xi}(\xi,t)}{-(n-2)\chi_\xi(\xi,t)}
	= \frac{\Phi''(z)}{-(n-2)\sqrt{t+t_0} \Phi'(z)} 
	\le \frac{c_2}{(n-2)c_1 (\xi+\xi_0)}.
  \eas
  Since $\xi+\xi_0>z_0\sqrt{t+t_0}$, (\ref{2.8}) implies that
  \bas
	\frac{\chi_{\xi\xi}(\xi,t)}{-(n-2)\chi_\xi(\xi,t)}
	< \frac{c_2}{(n-2)c_1 z_0 \sqrt{t+t_0}} < \frac{c_2}{(n-2) c_1 z_0 \sqrt{t_\star}} \le 1
  \eas
  holds at any such point, as claimed.
\qed
\begin{lem}\label{lem6}
  Let $D>0$ and $\gamma \in (0,1)$. Then there exists $t_\star>1$ such that for any choice of $t_0>t_\star$,
  $\xi_0>0$ and $A>0$, the function $\chip$ in {\rm (\ref{2.1})} satisfies
  \be{6.1}
	\qarab \chip \ge 0
	\qquad \mbox{for all $\xi>0$ and $t>0$,}
  \ee
  where $\qarab$ is the operator defined in {\rm (\ref{0chi})}.
\end{lem}
\proof
  We take $t_\star$ as given by Lemma~\ref{lem2} and assume that $t_0>t_\star$. Then writing
  $\chi:=\chip$ and using (\ref{2.2}) and (\ref{2.3}), we obtain
  \bas
	\qarab \chi &=& - e^{-2\xi} \Big\{ (D+\chi)[\chi_{\xi\xi} + (n-2)\chi_\xi] - \frac{n-2}{2} \chi_\xi^2 \Big\} \\
	&\ge& e^{-2\xi} \frac{n-2}{2} \chi_\xi^2 
	\ge 0
	\qquad \mbox{for all $\xi>0$ and $t>0$,}
  \eas
  because $D+\chi \ge 0$ according to the nonnegativity of $\chi$ asserted by Lemma~\ref{lem1}, and because $n\ge 3$.
\qed
The function we shall use as a supersolution near the origin (cf.~(\ref{phip}) below) 
will have a certain self--similar structure.
%
%
As a preparation, let us state the following.
\begin{lem}\label{lem4}
  Let $D>0$ and $\gamma>0$. For $\lambda:=\frac{1}{D}(\frac{\gamma}{2}+1)$, let $\rho$ denote the solution of
  \be{4.1}
	\left\{ \begin{array}{l}
	\rho''(\sigma)+ \frac{1}{\sigma} \rho'(\sigma) + \lambda \rho(\sigma)=0, \qquad \sigma>0, \\[1mm]
	\rho(0)=1, \quad \rho'(0)=0.
	\end{array} \right.
  \ee
  Then there exists $\sigma_0 \in (0,1)$ such that $\rho$ is positive and decreasing on $[0,\sigma_0]$.
\end{lem}
\proof
  Both statements are obvious from (\ref{4.1}).
\qed
In order to match inner and outer functions appropriately, we shall need a correcting factor which is time--dependent 
but approaches one in the large time limit.
\begin{lem}\label{lem3}
  Given $D>0$ and $\gamma \in (0,1)$, let $\Phi$, $\rho$ and $\sigma_0$ 
  be as in Lemma~\ref{lem1} and Lemma~\ref{lem4}, respectively.
  Then for $\xi_0>0$ and $t_0>\sigma_0^{-2}$, the function $\fp$ defined by
  \be{3.1}
	\fp(t):=\frac{\Phi\big(\xi_0(t+t_0)^{-\frac{1}{2}}\big)}{\rho\big((t+t_0)^{-\frac{1}{2}}\big)}, \qquad t\ge 0,
  \ee
  satisfies
  \be{3.2}
	\fp(t) \to 1 \qquad \mbox{as } t\to\infty.
  \ee
  Furthermore, for any $\xi_0>0$ there exists $C(\xi_0)>0$ such that whenever $t_0>1$, we have
  \be{3.3}
	\Big| (\fp)'(t) \Big| \le \frac{C(\xi_0)}{(t+t_0)^2}
	\qquad \mbox{for all } t>0.
  \ee
\end{lem}
\proof
  Since $\Phi(0)=\rho(0)=1$, (\ref{3.2}) is obvious.
  As for (\ref{3.3}), we for $t>0$ compute
  \bea{3.4}
	(\fp)'(t)
	=\frac{1}{2}(t+t_0)^{-\frac{3}{2}}\bigg(
&-&\xi_0
\frac{\Phi'\big(\xi_0(t+t_0)^{-\frac{1}{2}}\big)}{\rho\big((t+t_0)^{-\frac{1}{2}}\big)}\nn \\
		&+& \frac{\Phi\big(\xi_0(t+t_0)^{-\frac{1}{2}}\big) \rho'\big((t+t_0)^{-\frac{1}{2}}\big) }
			{\rho^2\big((t+t_0)^{-\frac{1}{2}}\big)}\bigg).
  \eea
 Since $\rho$ is positive on $[0,\sigma_0]$ and $\Phi'(0)=\rho'(0)=0$, we can pick 
  $c_1>0, c_2>0$ and $c_3>0$ such that
  \bas
	\rho(\sigma) \ge c_1 \qquad \mbox{for all } \sigma \in [0,\sigma_0]
  \eas
  as well as
  \bas
	|\Phi'(z)| \le c_2 z
	\quad \mbox{for all } z\in [0,\xi_0]
	\quad \mbox{and} \quad
	|\rho'(\sigma)| \le c_3 \sigma 
	\quad \mbox{for all } \sigma \in [0,\sigma_0].
  \eas
  We thereby obtain from (\ref{3.4}) that for any choice of $t_0>\sigma_0^{-2}$ one has
  \bas
	\Big| (\fp)'(t)\Big|
	\le
\bigg(\frac{c_2 \xi_0^2}{2c_1} + \frac{c_3}{2c_1^2}\bigg)(t+t_0)^{-2}
	\qquad \mbox{for all } t>0,
  \eas
  because $\Phi \le 1$ on $[0,\infty)$ by Lemma~\ref{lem1}. 
\qed
We can now introduce a family of functions one of which will serve as a supersolution in the 
region where $r<1$. To this end, for $D>0$ and $\gamma \in (0,1)$ we let $\rho, \sigma_0$
and $\fp$ as in Lemma~\ref{lem4} and Lemma~\ref{lem3}, and given 
$\xi_0>0$, $t_0>\sigma_0^{-2}$ and $A>0$ we define for $ r\in [0,1]$ and $t\ge 0$ the function
\be{phip}
	\phip(r,t):=A\fp(t) (t+t_0)^{-\frac{\gamma}{2}} \rho \big(r(t+t_0)^{-\frac{1}{2}}\big).
\ee
We then have the following.
\begin{lem}\label{lem5}
  Let $D>0$ and $\gamma \in (0,1)$, and let $\rho$ and $\sigma_0$ be as in Lemma~\ref{lem4}. 
  Then for each $\xi_0>0$ there exists $t^{\star} > \sigma_0^{-2}$ such that for any choice
  of $t_0>t^{\star}$ and any $A>0$, the function $\phip$ given by {\rm (\ref{phip})} satisfies
  \be{5.2}
	\parab \phip \ge 0
	\qquad \mbox{for all $r\in (0,1)$ and $t>0$.}
  \ee
\end{lem}
\proof
  Given $\xi_0>0$, we take $C(\xi_0)$ as provided by Lemma~\ref{lem3}, and claim that (\ref{5.2}) is valid whenever
  $t_0>t^{\star}$ and
  \be{5.22}
	t^{\star}>\max \Big\{ \frac{1}{\sigma_0^2}, \frac{C(\xi_0)}{\Phi(\xi_0)} \Big\},
  \ee
  where $\Phi$ is from Lemma~\ref{lem1}.

  To see this, we fix any such $t_0$ and, writing $\varphi=\phip$, $f=\fp$ and $\sigma=r(t+t_0)^{-1/2}$, compute
  \be{5.3}
	\varphi_r=Af(t) (t+t_0)^{-\frac{\gamma}{2}-\frac{1}{2}} \rho'(\sigma),\quad
	\varphi_{rr}=Af(t) (t+t_0)^{-\frac{\gamma}{2}-1} \rho''(\sigma)
  \ee
  and
  \bas
	\varphi_t = -\frac 12 Af(t) (t+t_0)^{-\frac{\gamma}{2}-1} (\sigma \rho'(\sigma) +
	\gamma \rho(\sigma))
	 - Af'(t) (t+t_0)^{-\frac{\gamma}{2}} \rho(\sigma)
  \eas
  for $r\in (0,1)$ and $t>0$.
  Now since $t_0>t^{\star}>\sigma_0^{-2}$, in the region where $r<1$ and $t>0$ we have
  $\sigma <t_0^{-1/2} < \sigma_0$, so that Lemma~\ref{lem4} guarantees that
  $\rho(\sigma)> 0$ and $\rho'(\sigma) \le 0$ and hence 
  $\rho''(\sigma) + \frac{1}{\sigma} \rho'(\sigma)=-\lambda \rho(\sigma) <0$.
 In particular, if we write  (\ref{0phi}) as
  \bas
	\parab \varphi = \varphi_t - (D+\varphi) \Big( \varphi_{rr} + \frac{n-1}{r} \varphi_r \Big)
	- r^2 \varphi_{rr} - r\varphi_r + \frac{n-2}{2} \varphi_r^2,
  \eas
  and use (\ref{5.3}), we obtain
  \bas
	-(D+\varphi) \Big(\varphi_{rr}&+&\frac{n-1}{r}\varphi_r\Big) \\
	&=& -(D+\varphi) Af(t) (t+t_0)^{-\frac{\gamma}{2}-1}
		\Big(\rho''(\sigma)+\frac{n-1}{\sigma} \rho'(\sigma) \Big) \\
	&=& - (D+\varphi) Af(t) (t+t_0)^{-\frac{\gamma}{2}-1}
		\Big(-\lambda \rho(\sigma) + \frac{n-2}{\sigma}\rho'(\sigma)\Big) \\
	&\ge& \lambda (D+\varphi) Af(t) (t+t_0)^{-\frac{\gamma}{2}-1} \rho(\sigma) \\
	&\ge& \lambda D Af(t) (t+t_0)^{-\frac{\gamma}{2}-1} \rho(\sigma),
  \eas
  because $n\ge 3$. Moreover,
  \bas
	-r^2 \varphi_{rr}-r\varphi_r
	&=& -Af(t) (t+t_0)^{-\frac{\gamma}{2}-1} r^2 \Big( \rho''(\sigma)+\frac{1}{\sigma} \rho'(\sigma)\Big) \\
	&=& \lambda A f(t) (t+t_0)^{-\frac{\gamma}{2}-1} r^2 
	\ge 0
  \eas
  for $r<1$ and $t>0$. Since $\frac{n-2}{2} \varphi_r^2 \ge 0$, we therefore have
  \bea{5.5}
	\parab \varphi &\ge& \varphi_t + \lambda DA f(t) (t+t_0)^{-\frac{\gamma}{2}-1} \rho(\sigma) \nn\\
	&=& Af(t) (t+t_0)^{-\frac{\gamma}{2}-1} 
		\Big\{ - \frac{\sigma}{2} \rho'(\sigma) - \frac{\gamma}{2} \rho(\sigma) 
		- \frac{f'(t)}{f(t)} (t+t_0) \rho(\sigma) \nn\\
	& & \hspace*{26mm}
	+ \lambda D \rho(\sigma) \Big\}
	\qquad \mbox{for $r\in (0,1)$ and } t>0.
  \eea

  Now the monotonicity properties of $\Phi$ and $\rho$ imply that, since $t_0>1$, we obtain
  \bas
	f(t) \ge \frac{1}{\rho(0)} \Phi\Big(\xi_0 t_0^{-\frac{1}{2}}\Big)
	\ge \Phi(\xi_0)
	\qquad \mbox{for all } t>0,
  \eas
  so that using (\ref{3.3}) we obtain
  \bas
	\Big| \frac{f'(t)}{f(t)} (t+t_0)\Big|
	\le \frac{C(\xi_0)}{\Phi(\xi_0) (t+t_0)}
	\le \frac{C(\xi_0)}{\Phi(\xi_0) t_0}
	\qquad \mbox{for all } t>0.
  \eas
  Thus, according to the fact that $t^{\star}>C(\xi_0)/\Phi(\xi_0)$ by (\ref{5.22}), we have
\bas
\Big|\frac{f'(t)}{f(t)} (t+t_0)\Big| \le 1 \qquad \mbox{for all } t>0.
\eas
 Hence, (\ref{5.5}) entails that
  \bas
	\parab \varphi &\ge& Af(t) (t+t_0)^{-\frac{\gamma}{2}-1}
	\bigg\{ - \frac{\sigma}{2} \rho'(\sigma) + \Big(\lambda D - \frac{\gamma}{2}-1\Big) \rho(\sigma) \bigg\} \\
	&=& -Af(t) (t+t_0)^{-\frac{\gamma}{2}-1} \frac{\sigma}{2} \rho'(\sigma) 
	\ge 0
	\qquad \mbox{for $r\in (0,1)$ and } t>0,
  \eas
  because of our choice of $\lambda$ in (\ref{4.1}) and, again, the monotonicity of $\rho$ on $(0,\sigma_0)$. 
  This completes the proof.
\qed
\begin{lem}\label{lem7}
  Let $D>0$ and $\gamma \in (0,1)$. Then with $\sigma_0$ as in Lemma~\ref{lem4}, there exist $\xi_0>0$ and 
  $t_0>\sigma_0^{-2}$ such that for any $A>0$, the function $\ophia$ defined by
  \be{7.1}
	\ophia(r,t):=\left\{ \begin{array}{ll}
	\phip(r,t), \qquad & r \in [0,1], \, t\ge 0, \\[1mm]
	\chip(\ln r,t), \qquad & r>1, \, t\ge 0,
	\end{array} \right.
  \ee
  is continuous in $[0,\infty)^2$ and satisfies
  \be{7.2}
	\parab \ophia \ge 0
	\qquad \mbox{for all $r\in (0,\infty) \setminus \{1\}$ and } t>0,
  \ee
  where $\parab$ is as in {\rm (\ref{0phi})}, and such that
  \be{7.3}
	\liminf_{r\nearrow 1} \ophia_r(r,t) > \limsup_{r\searrow 1} \ophia_r(r,t)
	\qquad \mbox{for all } t>0.
  \ee
\end{lem}
\proof
  Given $D>0$ and $\gamma \in (0,1)$, we let $\rho$ and $\Phi$ be as defined by (\ref{4.1}) and (\ref{1.1}).
  Then since $\rho'(0)=\Phi'(0)=0$ and $\Phi''(0)=-\frac{\gamma}{2}<0$, we can find $c_1>0$ and $c_2>0$ fulfilling
  \be{7.4}
	\rho'(\sigma) \ge - c_1 \sigma 
	\qquad \mbox{for all } \sigma \in (0,\sigma_0)
  \ee
  and
  \be{7.5}
	\Phi'(z) \le - c_2 z
	\qquad \mbox{for all } z\in (0,1).
  \ee
  We now first fix $\xi_0>0$ large such that
  \be{7.6}
	\xi_0>\frac{c_1}{c_2 \rho(\sigma_0)}
  \ee
  and then take $t_\star$ and $t^{\star}$ as provided by Lemma~\ref{lem6} and Lemma~\ref{lem5}, respectively,
  when applied to this particular choice of $\xi_0$. 
  We finally pick some $t_0>\sigma_0^{-2}$ satisfying
  \be{7.7}
	t_0>\max \big\{ t_\star, t^{\star}, \xi_0^2 \big\}
  \ee
  and claim that these choices ensure that $\ophia$ is continuous and that (\ref{7.2}) and (\ref{7.3}) are valid
  whenever $A>0$.

  In fact, (\ref{7.2}) is an immediate consequence of Lemma~\ref{lem6} 
and Lemma~\ref{lem5}, while the continuity of
  $\ophia$ directly results from the definitions of $\phip$, $\chip$ and the function $\fp$ introduced in Lemma~\ref{lem3}. 
  To  verify (\ref{7.3}), we recall (\ref{phip}) and (\ref{2.1}) in computing
  \bea{7.8}
	I_1(t) &:=& \liminf_{r\nearrow 1} \ophia_r(r,t) 
	= \phip_r(1,t) \nn\\
	&=& A \fp(t) (t+t_0)^{-\frac{\gamma}{2}-\frac{1}{2}} \rho'\left((t+t_0)^{-\frac{1}{2}}\right),
	\qquad t>0,
  \eea
  and
  \bea{7.9}
	I_2(t) &:=& \limsup_{r\searrow 1} \ophia_r(r,t) 
	= \chip_\xi(0,t) \nn\\
	&=& A (t+t_0)^{-\frac{\gamma}{2}-\frac{1}{2}} \Phi' \left( \xi_0 (t+t_0)^{-\frac{1}{2}}\right),
	\qquad t>0.
  \eea
  Here we note that by (\ref{3.1}) and the monotonicity of $\Phi$ and $\rho$,
  \be{7.10}
	\fp(t) \le \Phi(0)\rho^{-1}\left(t_0^{-\frac{1}{2}}\right)
	=\rho^{-1}\left(t_0^{-\frac{1}{2}}\right)
	\le \rho^{-1}(\sigma_0)
	\qquad \mbox{for all } t>0,
  \ee
  because $t_0>\sigma_0^{-2}$. Furthermore, (\ref{7.4}) and (\ref{7.5}) assert that
  \be{7.11}
	\rho'\left((t+t_0)^{-\frac{1}{2}}\right) \ge - c_1(t+t_0)^{-\frac{1}{2}}
	\qquad \mbox{for all } t>0
  \ee
  and 
  \be{7.12}
	\Phi' \left( \xi_0 (t+t_0)^{-\frac{1}{2}}\right) \le - c_2 \xi_0 (t+t_0)^{-\frac{1}{2}}
	\qquad \mbox{for all } t>0,
  \ee
  again since $(t+t_0)^{-1/2}<\sigma_0$, and since $\xi_0(t+t_0)^{-1/2}<1$  due to  
  (\ref{7.7}). 
  Using (\ref{7.10})--(\ref{7.12}), we obtain  from (\ref{7.8}) and (\ref{7.9}) that
  \bas
	I_1(t)-I_2(t)
	\ge A (t+t_0)^{-\frac{\gamma}{2}-1}
		\big( - c_1 \rho^{-1}(\sigma_0)
		+ c_2 \xi_0 \big)
	\quad \mbox{for all } t>0,
  \eas
  so that our requirement (\ref{7.6}) guarantees that (\ref{7.3}) holds.
\qed
\begin{lem}\label{lem8}
 Let $D>0$. Assume that $\varphi_0$ is continuous and nonnegative on $[0,\infty)$ and there exist 
  $\gamma \in (0,1)$ and $B>0$  such that
  \be{8.1}
	\varphi_0(r) \le B \ln^{-\gamma} r
	\qquad \mbox{for all } r>2.
  \ee
  Then there exists $C>0$ such that the solution $\varphi$ of {\rm (\ref{0phi})} with $\varphi(\cdot,0)=\varphi_0$ satisfies
  \be{8.2}
	\varphi(r,t) \le C (t+1)^{-\frac{\gamma}{2}}
	\qquad \mbox{for all $r\ge 0$ and } t\ge 0.
 \ee
\end{lem}
\proof
  Given $D>0$ and $\gamma \in (0,1)$, we fix $\sigma_0 \in (0,1)$, $\xi_0>0$ and $t_0>\sigma_0^{-2}$
  as in Lemma~\ref{lem4} and Lemma~\ref{lem7} and take $f=\fp$ from Lemma~\ref{lem3}. 
  In order to define, with some specific $A>0$, a supersolution of the form (\ref{7.1}) which initially dominates $\varphi$,
  we set $z_0:=(\ln 2 + \xi_0)t_0^{-1/2}$ and then obtain from Lemma~\ref{lem1} that for some $c_1>0$,
  the function $\Phi$ in (\ref{1.1}) satisfies
  \be{8.3}
	\Phi(z) \ge c_1 z^{-\gamma}
	\qquad \mbox{for all } z\ge z_0.
  \ee
  Moreover, since $\varphi_0$ is bounded we can pick $c_2>0$ such that
  \be{8.4}
	\varphi_0(r) \le c_2
	\qquad \mbox{for all } r\in [0,2].
  \ee
  We now fix any $A>0$ fulfilling
  \be{8.5}
	A>\max \left\{ \frac{c_2 t_0^{\gamma/2}}{f(0) \rho(t_0^{-1/2})} \, , \,
	\frac{c_2 t_0^{\gamma/2}}{\Phi(z_0)} \, , \, \frac{B}{c_1}  \left(1+\frac{\xi_0}{\ln 2}\right)^{\gamma}
	\right\}
  \ee
  and claim that then the function $\ophia$ in (\ref{7.1}) has the property 
  \be{8.6}
	\ophia(r,0) >\varphi_0(r)
	\qquad \mbox{for all } r\ge 0.
  \ee
  To prove this, we first observe that for small $r$, by (\ref{8.4}) and (\ref{8.5}) it holds that
  \bas
	\frac{\ophia(r,0)}{\varphi_0(r)}
	&\ge& \frac{\ophia(r,0)}{c_2}
	= \frac{1}{c_2}Af(0) t_0^{-\frac{\gamma}{2}} \rho \big(rt_0^{-\frac{1}{2}}\big) \\
	&\ge& \frac{1}{c_2}Af(0) t_0^{-\frac{\gamma}{2}} \rho \big(t_0^{-\frac{1}{2}}\big)
	>1,
	\qquad r\in [0,1],
  \eas
  because $\rho' \le 0$ on $(0,\sigma_0)$ and $t_0^{-1/2} < \sigma_0$. Similarly, in the intermediate region
  where $1<r \le 2$, (\ref{8.4}), (\ref{8.5}) and the monotonicity of $\Phi$ yield
  \bas
	\frac{\ophia(r,0)}{\varphi_0(r)}
	\ge \frac{1}{c_2}At_0^{-\frac{\gamma}{2}} \Phi\Big((\ln 2 + \xi_0)\big(t_0^{-\frac{1}{2}}\big)\Big) > 1,
	\qquad r \in (1,2].
  \eas
  Finally, for large $r$ we apply (\ref{8.3}) to estimate
  \bas
	\ophia(r,0)
	&=& A t_0^{-\frac{\gamma}{2}} \Phi\Big((\ln r + \xi_0)\big(t_0^{-\frac{1}{2}}\big)\Big)
	\ge  c_1 A (\ln r + \xi_0)^{-\gamma}\\
	&\ge& c_1 A \Big(1+\frac{\xi_0}{\ln 2} \Big)^{-\gamma} (\ln r)^{-\gamma},
	\qquad r>2,
  \eas
  because $\ln r + \xi_0 \le \ln r + \xi_0 \ln r/\ln 2$ for such $r$.
  Along with (\ref{8.5}) and (\ref{8.1}), this guarantees that also
  \bas
	\ophia(r,0) > \varphi_0(r), \qquad r>2.
  \eas
  Having thus found that (\ref{8.6}) is true, we may invoke Lemma~\ref{lem7} combined with the comparison principle
  to infer that $\varphi \le \ophia$ in $[0,\infty)^2$. In particular, since $\rho \le 1, \Phi \le 1$ and 
  $f \le  \rho^{-1}(t_0^{-1/2})$
 by monotonicity, this means that
  \bas
	\varphi(r,t) \le Af(t) (t+t_0)^{-\frac{\gamma}{2}}
	\le A \rho^{-1} \big(t_0^{-\frac{1}{2}}\big) (t+t_0)^{-\frac{\gamma}{2}},
	\qquad r\in [0,1],\quad t\ge 0,
  \eas
  as well as
  \bas
	\varphi(r,t) \le A(t+t_0)^{-\frac{\gamma}{2}}
	\qquad \mbox{for all $r>1$ and } t\ge 0,
  \eas
  from which (\ref{8.2}) clearly follows.
\qed
\medskip
\noindent
{\it Proof of Theorem~\ref{thm1}} (i)
If we choose $\varphi_0$ satisfying (\ref{8.1}) such that $\psi_0(x)\le\varphi_0(|x|)$
for $x\in\R^n$, we obtain by comparison that 
\bas
\left(|x|^2 + D + \varphi(|x|,t)\right)^{-\frac{n-2}{2}}\le v(x,t)\le V_D(x),
\qquad x\in\R^n,\quad t\ge 0.
\eas
Lemma~\ref{lem8} and the Mean Value Theorem yield then the result.
\qed
\mysection{Upper bound. Proof of Theorem \ref{thm1} (ii)}\label{Sec3}
\begin{lem}\label{lem100}
  For $\gamma>0$, let
  \be{100.1}
	\hphi(z):= \Big(1+\frac{z^2}{4} \Big)^{-\frac{\gamma}{2}}, \qquad z\ge 0.
  \ee
  Then 
  \be{100.2}
	\hphi''(z) + \frac{z}{2} \hphi'(z) + \frac{\gamma}{2} \hphi(z) 
	\ge \frac{\gamma}{4} \Big(1+\frac{z^2}{4} \Big)^{-\frac{\gamma}{2}-1}
	\qquad \mbox{for all } z>0.
  \ee
  Moreover,
  \be{100.3}
	\hphi'(z)= -\frac{\gamma}{4} z \Big(1+\frac{z^2}{4} \Big)^{-\frac{\gamma}{2}-1}
	\qquad \mbox{for all } z>0
  \ee
  and 
  \be{100.4}
	|\hphi''(z)| \le \frac{\gamma(\gamma+1)}{4} \Big(1+\frac{z^2}{4} \Big)^{-\frac{\gamma}{2}-1}
	\qquad \mbox{for all } z>0.
  \ee
\end{lem}
\proof
  By straightforward computation, we find (\ref{100.3}) as well as
  \be{100.6}
	\hphi''(z) = \frac{\gamma(\gamma+2)}{16} z^2 \Big(1+\frac{z^2}{4} \Big)^{-\frac{\gamma}{2}-2}
	-\frac{\gamma}{4} \Big(1+\frac{z^2}{4} \Big)^{-\frac{\gamma}{2}-1},
  \ee
  and thus obtain
  \bas
	I := \hphi''(z) + \frac{z}{2} \hphi'(z) + \frac{\gamma}{2} \hphi(z) 
	= \Big(1+\frac{z^2}{4} \Big)^{-\frac{\gamma}{2}-2} \bigg\{
	\frac{\gamma}{4} + \frac{\gamma(\gamma+3)}{16} z^2 \bigg\}
  \eas
  for $z>0$. 
This implies that
  \bas
	I \ge \frac{\gamma}{4} \Big(1+\frac{z^2}{4} \Big)^{-\frac{\gamma}{2}-1}
	\qquad \mbox{for all } z>0
  \eas
  and thereby proves (\ref{100.2}). 
  Moreover, by (\ref{100.6}) we infer that
  \bas
	|\hphi''(z)| 
	&=& \Big(1+\frac{z^2}{4} \Big)^{-\frac{\gamma}{2}-2} 
		\Big| \frac{\gamma(\gamma+1)}{16} z^2 + \frac{\gamma}{4} \Big| \\
	&\le& \Big(1+\frac{z^2}{4} \Big)^{-\frac{\gamma}{2}-2}
	\Big| \frac{\gamma(\gamma+1)}{16} z^2 + \frac{\gamma(\gamma+1)}{4} \Big| \\
	&=& \frac{\gamma(\gamma+1)}{4} \Big(1+\frac{z^2}{4} \Big)^{-\frac{\gamma}{2}-1} 
		\qquad \mbox{for all } z>0.
  \eas
\qed
\begin{lem}\label{lem101}
  Let $D>0$ and $\gamma>0$. Then there exists $\xi_0>1$ such that for any choice of $a\in (0,1)$, the function
  $\hchia$ given by
  \be{101.1}
	\hchia(\xi,t):=a(t+1)^{-\frac{\gamma}{2}} \hphi \Big(\frac{\xi-\xi_0}{\sqrt{t+1}} \Big),
	\qquad \xi\ge \xi_0, \ t\ge 0,
  \ee
  satisfies
  \be{101.2}
	\qarab \hchia \le 0
	\qquad \mbox{for all $\xi>\xi_0$ and } t>0,
  \ee
  where $\qarab$ is as defined in {\rm (\ref{0chi})}.
\end{lem}
\proof
  We abbreviate 
  \bas
	z:=\frac{\xi-\xi_0}{\sqrt{t+1}}
  \eas
  and calculate
  \bas
	\hchia_\xi = a (t+1)^{-\frac{\gamma}{2}-\frac{1}{2}} \hphi'(z), \qquad
	\hchia_{\xi\xi} = a (t+1)^{-\frac{\gamma}{2}-1} \hphi''(z)
  \eas
  and
  \bas
	\hchia_t = -\frac{1}{2} a (t+1)^{-\frac{\gamma}{2}-1} z\hphi'(z)
	-\frac{\gamma}{2} a (t+1)^{-\frac{\gamma}{2}-1} \hphi(z).
  \eas
  Therefore,
  \bas
	\qarab \hchia &=& a(t+1)^{-\frac{\gamma}{2}-1} \bigg\{
	-\frac{z}{2} \hphi'(z) - \frac{\gamma}{2} \hphi(z) - \hphi''(z) \\
	& & \hspace*{10mm}
	- e^{-2\xi} \Big[ D+a(t+1)^{-\frac{\gamma}{2}} \hphi(z) \Big] 
		\Big[ \hphi''(z) + (n-2) \sqrt{t+1} \hphi'(z) \Big] \\
	& & \hspace*{10mm}
	+ \frac{n-2}{2} e^{-2\xi} (t+1)^{-\frac{\gamma}{2}} \hphi'^2(z) \bigg\}
	\qquad \mbox{for $\xi>\xi_0$ and } t>0.
  \eas
  Here we recall (\ref{100.1}) and our assumption $a<1$ in estimating
  \bas
	\Big| D + a(t+1)^{-\frac{\gamma}{2}} \hphi(z) \Big| \le D+1
  \eas
  and use (\ref{100.2}) to see that
  \bas
	-\frac{z}{2} \hphi'(z) - \frac{\gamma}{2} \hphi(z) - \hphi''(z) 
	\le -\frac{\gamma}{4} \Big(1+\frac{z^2}{4} \Big)^{-\frac{\gamma}{2}-1}
  \eas
  at any point $(\xi,t) \in (\xi_0,\infty) \times (0,\infty)$. Thus,
  \bas
	\frac{\qarab \hchia}{a(t+1)^{-\frac{\gamma}{2}-1}}
	&\le& -\frac{\gamma}{4} \Big(1+\frac{z^2}{4}\Big)^{-\frac{\gamma}{2}-1} 
	 + (D+1) e^{-2\xi} |\hphi''(z)|\\
&& + (n-2) e^{-2\xi} \sqrt{t+1} |\hphi'(z)| 
	 + \frac{n-2}{2} e^{-2\xi} \hphi'^2(z) \\
	&=:& -I_1+I_2+I_3+I_4
	\qquad \mbox{for $\xi>\xi_0$ and } t>0,
  \eas
  and we claim that this implies (\ref{101.2}) if we pick $\xi_0>1$ large enough such that
  \be{101.3}
	(\gamma+1)(D+1) e^{-2\xi} < \frac{1}{3}
	\qquad \mbox{for all } \xi>\xi_0
  \ee
  and
  \be{101.4}
	(n-2)(D+1) (\xi-\xi_0) e^{-2\xi} < \frac{1}{3}
	\qquad \mbox{for all } \xi>\xi_0
  \ee
  as well as
  \be{101.5}
	\frac{(n-2)\gamma}{2} e^{-2\xi} < \frac{1}{3}
	\qquad \mbox{for all } \xi>\xi_0.
  \ee
  Indeed, in conjunction with (\ref{100.4}), (\ref{101.3}) implies that
  \bea{101.6}
	\frac{I_2}{I_1}
	&=& \frac{4(D+1)}{\gamma} e^{-2\xi} \Big(1+\frac{z^2}{4} \Big)^{\frac{\gamma}{2}+1} |\hphi''(z)| \nn\\
	&\le& \frac{4(D+1)}{\gamma} e^{-2\xi} \frac{\gamma(\gamma+1)}{4} 
	= (\gamma+1)(D+1) e^{-2\xi} \nn\\
	&<& \frac{1}{3}
	\qquad \mbox{for all $\xi>\xi_0$ and } t>0.
  \eea
  Since $\sqrt{t+1}=(\xi-\xi_0)/z$, (\ref{100.3}) and (\ref{101.4}) next guarantee that
  \bea{101.7}
	\frac{I_3}{I_1}
	&=& \frac{4(n-2)(D+1)}{\gamma}  e^{-2\xi} \sqrt{t+1} \Big(1+\frac{z^2}{4} \Big)^{\frac{\gamma}{2}+1}
		|\hphi'(z)| \nn\\
	&=& \frac{4(n-2)(D+1)}{\gamma}  (\xi-\xi_0) e^{-2\xi} 
		\frac{(1+\frac{z^2}{4})^{\frac{\gamma}{2}+1}}{z} |\hphi'(z)| \nn\\
	&=& (n-2)(D+1)  (\xi-\xi_0) e^{-2\xi} \nn\\
	&<& \frac{1}{3}
	\qquad \mbox{for all $\xi>\xi_0$ and } t>0.
  \eea
  Finally, again by (\ref{100.3}),
  \bas
	\frac{I_4}{I_1}
	&=& \frac{2(n-2)}{\gamma}  e^{-2\xi}  (t+1)^{\frac{\gamma}{2}+1}  \hphi'^2(z) \\
	&=& \frac{(n-2)\gamma}{8}  e^{-2\xi}  z^2 \Big(1+\frac{z^2}{4} \Big)^{-\frac{\gamma}{2}-1},
  \eas
  so that, since clearly $z^2 (1+\frac{z^2}{4})^{-\frac{\gamma}{2}-1} \le 4$, from (\ref{101.5}) we infer that
  \bas
	\frac{I_4}{I_1} \le \frac{(n-2)\gamma}{2}  e^{-2\xi} < \frac{1}{3}
	\qquad \mbox{for all $\xi>\xi_0$ and } t>0.
  \eas
  Combined with (\ref{101.6}) and (\ref{101.7}), this establishes (\ref{101.2}).
\qed
%
%
%
%
  In view of the explicit definition (\ref{100.1}) of $\hphi$, the above function $\hchia$ can alternatively be 
  written in the fully explicit form
  \bas
	\hchia(\xi,t)=a \Big(t+1+\frac{(\xi-\xi_0)^2}{4} \Big)^{-\frac{\gamma}{2}},
	\qquad \xi\ge \xi_0, \ t\ge 0.\abs
  \eas
\begin{lem}\label{lem102}
  Let $D>0$ and $\gamma>0$. Then there exists $r_0>e$ such that for all $a\in (0,1)$,
  \be{102.1}
	\uphia(r,t):= \left\{ \begin{array}{ll}
	a(t+1)^{-\frac{\gamma}{2}}, \qquad & r \in [0,r_0], \ t\ge 0, \\[1mm]
	\hchia(\ln r,t), \qquad & r>r_0, \ t\ge 0,
	\end{array} \right.
  \ee
  defines a continuous function $\uphia$ on $[0,\infty)^2$ such that also $\uphia_r$ is continuous
  on $[0,\infty)^2$, and such that
  \be{102.2}
	\parab \uphia \le 0
	\qquad \mbox{for all $r\in (0,\infty) \setminus \{r_0\}$ and } t>0.
  \ee
  Here $\hchia$ is as defined in Lemma~\ref{lem101} with $\xi_0:=\ln r_0$, and $\parab$ is as in {\rm (\ref{0phi})}.
\end{lem}
\proof
  With $\xi_0>1$ as provided by Lemma~\ref{lem101}, we let $r_0:=e^{\xi_0}>e$ and thereupon obtain that
  (\ref{101.2}) precisely yields $\parab \uphia \le 0$ for $r>r_0$ and $t>0$. To see the same for small $r$,
  we only need to note that clearly
  \be{102.4}
	\uphia_r=\uphia_{rr} \equiv 0 
	\qquad \mbox{for $r<r_0$ and } t>0,
  \ee
  so that
  \bas
	\parab \uphia=\uphia_t = -\frac{a\gamma}{2} (t+1)^{-\frac{\gamma}{2}-1} < 0
	\qquad \mbox{for $r<r_0$ and } t>0.
  \eas
  Having thus established (\ref{102.2}), we are left with proving the continuity of $\uphia_r$.
  In view of (\ref{102.4}), however, this immediately follows from the observation that
  \bas
	\lim_{r\searrow r_0} \uphia_r(r,t)
	= \frac{1}{r_0} \hchia(\xi_0,t)
	= \frac{a}{r_0} (t+1)^{-\frac{\gamma}{2}} \hphi'(0) = 0
	\qquad \mbox{for all } t>0,
  \eas
  whereby the proof is completed.
\qed
\begin{lem}\label{lem103}
  Let $D>0$. Suppose that $\varphi_0 \in C^0([0,\infty))$ is positive and such that
  \be{103.1}
	\varphi_0(r) \ge b \ln^{-\gamma} r
	\qquad \mbox{for all } r>2
  \ee
 with some positive constants $b$ and $\gamma$. Then there exists $c>0$ such that the solution $\varphi$
  of {\rm (\ref{0phi})} fulfilling $\varphi(\cdot,0)=\varphi_0$ satisfies
  \be{103.2}
	\varphi(0,t) \ge c (t+1)^{-\frac{\gamma}{2}}
	\qquad \mbox{for all } t>0.
  \ee
\end{lem}
\proof
  We let $r_0>e$ be as given by Lemma~\ref{lem102}. Then since $\varphi_0$ is continuous and positive, we can find 
  $c_1>0$ such that
  \be{103.3}
	\varphi_0(r) \ge c_1
	\qquad \mbox{for all } r \in [0,r_0],
  \ee
  and fix $a \in (0,1)$ small enough fulfilling
  \be{103.4}
	a<\min \{c_1 \, , \, b c_2^\frac{\gamma}{2} \},
  \ee
  where  
  \bas
	c_2:=\min \Big\{ \frac{1}{16}, \frac{1}{4\xi_0^2} \Big\}
  \eas
  with $\xi_0:=\ln r_0>1$.
  We claim that this choice ensures that with $\uphia$ defined by (\ref{102.1}) we have
  \be{103.6}
	\varphi_0(r) \ge \uphia(r,0) 
	\qquad \mbox{for all } r\ge 0.
  \ee
  In fact, if $r$ is small then by (\ref{103.3}) and (\ref{103.4}),
  \bas
	\varphi_0(r) \ge c_1 > a =\uphia(r,0)
	\qquad \mbox{for all } r\in  [0,r_0].
  \eas
  In order to show (\ref{103.6}) for large $r$, we observe that by (\ref{102.1}), (\ref{101.1}) and (\ref{100.1}),
  \bas
	\uphia(r,0)
	= a\hphi(\ln r -\xi_0) 
	= a \Big(1+ \frac{(\ln r - \xi_0)^2}{4} \Big)^{-\frac{\gamma}{2}},
	\qquad r>r_0,
  \eas
  because $r_0>e>2$. Here we estimate
  \bas
	1+\frac{(\ln r - \xi_0)^2}{4} 
	\ge \frac{(\ln r-\xi_0)^2}{4} 
	\ge \frac{(\ln r)^2}{16}
	\qquad \mbox{if } \ln r \ge 2\xi_0
  \eas
  and 
  \bas
	1+\frac{(\ln r - \xi_0)^2}{4} 
	\ge 1
	\ge \Big( \frac{\ln r}{2\xi_0} \Big)^2
	\qquad \mbox{if } \ln r < 2\xi_0,
  \eas
  whence by definition of $c_2$ it follows that
  \bas
	\uphia(r,0) \le a \Big( c_2 (\ln r)^2 \Big)^{-\frac{\gamma}{2}}
	< b (\ln r)^{-\gamma}
	\le \varphi_0(r)
	\qquad \mbox{for all } r>2.
  \eas
  We have thereby verified (\ref{103.6}), which in turn, on an application of the comparison principle, entails that
  $\varphi \ge \uphia$ in $[0,\infty)^2$.
  Evaluated at $r=0$, this gives as a particular consequence that
  \bas
	\varphi(0,t) \ge \uphia(0,t) = a(t+1)^{-\frac{\gamma}{2}}
	\qquad \mbox{for all } t\ge 0
  \eas
  and hence proves (\ref{103.2}).
\qed
\medskip
\noindent
{\it Proof of Theorem~\ref{thm1}} (ii)
We choose $\varphi_0$ satisfying (\ref{103.1}) such that $\psi_0(x)\ge\varphi_0(|x|)$
for $x\in\R^n$. Then we obtain by comparison that 
\bas
\left(|x|^2 + D + \varphi(|x|,t)\right)^{-\frac{n-2}{2}}\ge v(x,t),
\qquad x\in\R^n,\quad t\ge 0.
\eas
Lemma~\ref{lem103} and the Mean Value Theorem yield then the result.
\qed
\mysection{Universal upper bound. Proof of Theorem \ref{thm2}}\label{Sec4}
\begin{lem}\label{lem201}
  Let $\xi_1\in\R$, and suppose that $\alpha$ and $\beta$ are smooth functions on $(\xi_1,\infty) \times (0,\infty)$,
  for which there exist $k>0$ and $K>0$ such that
  \bas
	k \le \alpha(\xi,t) \le K \quad \mbox{and} \quad |\beta(\xi,t)| \le K
	\qquad \mbox{for all $\xi>\xi_1$ and } t>0.
  \eas
  Then for any nonnegative solution 
$$0 \not\equiv w \in C^{2,1}((\xi_1,\infty) \times (0,\infty)) \cap
  C^0([\xi_1,\infty) \times [0,\infty))$$ of 
  \bas
	w_t=\alpha(\xi,t) w_{\xi\xi} + \beta(\xi,t) w_\xi,
	\qquad \xi>\xi_1, \, t>0,
  \eas
  one can find $c>0$ such that
  \bas
	\sup_{\xi>\xi_1} w(\xi,t) \ge c(t+1)^{-\frac{1}{2}}
	\qquad \mbox{for all } t>0.
  \eas
\end{lem}
\proof
 This lower bound follows from \cite{A}, for example.
\qed
\begin{lem}\label{lem200}
  Let $D>0$ and assume that $\varphi_0$ is continuous and nonnegative on $[0,\infty)$, $\varphi_0 \not\equiv 0$.
  Then there exists $c>0$ such that the solution $\varphi$ of {\rm (\ref{0phi})} with $\varphi(\cdot,0)=\varphi_0$ satisfies
  \be{200.1}
	\sup_{r>0} \varphi(r,t) \ge c(t+1)^{-\frac{1}{2}}
	\qquad \mbox{for all } t>0.
  \ee
\end{lem}
\proof
  Passing to a suitable minorant of $\varphi_0$ if necessary, in view of the comparison principle we may assume that for
  some $r_0>0$ we have $0 \not\equiv \varphi_0 \in C_0^\infty((r_0,\infty))$ with $0\le \varphi_0 \le 1$. 
  Now conveniently rewritten in terms of $\chi(\xi,t)=\varphi(r,t)$, $\xi=\ln r$, (\ref{0phi}) becomes
  (cf.~also (\ref{0chi}))
  \bas
	\chi_t &=& \chi_{\xi\xi}
	+ e^{-2\xi} \Big\{ (D+\chi) [\chi_{\xi\xi} + (n-2)\chi_\xi] - \frac{n-2}{2} \chi_\xi^2 \Big\} \\
	&=& \Big[1+(D+\chi) e^{-2\xi} \Big] \chi_{\xi\xi}
	+ e^{-2\xi} \Big[ (n-2)(D+\chi) - \frac{n-2}{2} \chi_\xi \Big] \chi_\xi \\
	&=:& \alpha(\xi,t) \chi_{\xi\xi} + \beta(\xi,t) \chi_\xi,
	\qquad \xi\in\R, \ t>0.
  \eas
Let us next choose $\xi_0\in\R$ such that $\xi_0<\ln r_0-2$.
  Then, since $0\le \chi \le 1$ in $\R \times (0,\infty)$, we have
  \bas
	1 \le \alpha(\xi,t) \le 1+ (D+1) e^{-2\xi_0}
	\qquad \mbox{for all $\xi>\xi_0$ and } t>0.
  \eas
  Therefore, due to the fact that $\varphi_0$ is smooth with compact support, interior parabolic Schauder estimates
  (\cite{LSU}) provide $c_1>0$ such that
  \bas
	|\chi_\xi(\xi,t)| \le c_1
	\qquad \mbox{for all $\xi>\xi_1$ and } t>0,
  \eas
  so that 
  \bas
	|\beta(\xi,t)| \le e^{-2(\xi_0+1)} \Big[ (n-2)(D+1) + \frac{n-2}{2} c_1 \Big]
	\quad \mbox{for all $\xi>\xi_1$ and } t>0.
  \eas
  Since we already know that $\chi \ge 0$ and that $\chi(\cdot,0) \not\equiv 0$ in $(\xi_0+1,\infty)$ according to our
  choices of $r_0$ and $\xi_0$, we may now invoke Lemma~\ref{lem201} to conclude that there exists $c_2>0$ such that
  \bas
	\sup_{\xi>\xi_0+1} \chi(\xi,t) \ge c_2 (t+1)^{-\frac{1}{2}}
	\qquad \mbox{for all } t>0.
  \eas
  Restated using the variable $\varphi$, this immediately yields (\ref{200.1}).
\qed
\medskip
\noindent
{\it Proof of Theorem~\ref{thm2}}
We can assume without loss of generality that $\psi_0(0)>0$.
We choose $\varphi_0$ such that $\psi_0(x)\ge\varphi_0(|x|)$ for $x\in\R^n$
and $\varphi_0\not\equiv 0$ is nonincreasing.
We then obtain by comparison that 
\bas
\left(|x|^2 + D + \varphi(|x|,t)\right)^{-\frac{n-2}{2}}\ge v(x,t),
\qquad x\in\R^n,\quad t\ge 0.
\eas
Since $\sup_{r>0} \varphi(r,t)=\varphi(0,t)$,
the result follows from Lemma~\ref{lem200} and the Mean Value Theorem.
\qed

\bigskip
\noindent
{\bf Acknowledgment.} \
The first author was supported in part by the Slovak Research and Development Agency under the contract
No. APVV-0134-10 and by the VEGA grant 1/0711/12. This work was finished 
while the first author was visiting
the Interactive Research Center of Science, Tokyo Institute of Technology. He is grateful for the warm hospitality.
The second author gratefully acknowledges the support of the Royal Society,
of the ESF network HCAA and Wolfson Foundation.

\end{document}